\documentclass{amsart}
\usepackage{amssymb,amscd,amsmath,eucal}
\usepackage[matrix,arrow,curve,frame]{xy}    
\xymatrixcolsep{1.9pc}                          
\xymatrixrowsep{1.9pc}
\newdir{ >}{{}*!/-8pt/\dir{>}}                  

\makeatletter
\renewcommand{\phi}{\varphi}
\renewcommand{\geq}{\geqslant}
\renewcommand{\leq}{\leqslant}
\renewcommand{\epsilon}{\varepsilon}

\renewcommand{\kappa}{\varkappa}

\DeclareMathOperator{\spec}{\sf Spec} \DeclareMathOperator{\ann}{\sf
ann} \DeclareMathOperator{\thick}{thick}
\DeclareMathOperator{\serre}{tor} \DeclareMathOperator{\inj}{\sf
Inj} \DeclareMathOperator{\Hom}{Hom} \DeclareMathOperator{\supp}{\sf
supp} \DeclareMathOperator{\open}{{open}}

\DeclareMathOperator{\injfg}{\sf Inj_{\rm fg}}
\DeclareMathOperator{\injzg}{\sf Inj_{\rm zg}}

\DeclareMathOperator{\End}{End}

 \DeclareMathOperator{\Gr}{Gr}
 \DeclareMathOperator{\gr}{gr}
 
 \DeclareMathOperator{\Qcoh}{Qcoh}
 \DeclareMathOperator{\kr}{Ker}
 
 \DeclareMathOperator{\im}{Im}
 \DeclareMathOperator{\QGr}{QGr}
 \DeclareMathOperator{\Tors}{Tors}
 
\DeclareMathOperator{\Proj}{\sf Proj}
\DeclareMathOperator{\perf}{\cc D_{per}}
 
\DeclareMathOperator{\Mod}{Mod}

\newcommand {\lp}{\varinjlim}

\newcommand{\lra}[1]{\bl{#1}\longrightarrow\relax}
\newcommand{\bl}[1]{\buildrel #1\over}
\newcommand{\cc}{\mathcal}
\newcommand{\ps}{\oplus}
\newcommand{\ff}{\mathfrak}

\newcommand{\wt}{\widetilde}
\newcommand{\ifff}{if and only if }

\newcommand{\bb}{\mathbb}

\newcommand{\Rfp}{\Mod R}

\newtheorem{thm}{Theorem}[section]
\newtheorem{prop}[thm]{Proposition}
\newtheorem{cor}[thm]{Corollary}
\newtheorem{lem}[thm]{Lemma}

\newtheorem*{theo}{Theorem}

\newtheorem{rem}[thm]{Remark}

\newtheorem*{defs}{Definition}

\numberwithin{equation}{section}

\makeatother
\begin{document}

\footskip30pt


\title{Torsion classes of finite type and spectra}
\author{Grigory Garkusha}
\address{School of Mathematics, University of Manchester, Oxford Road, Manchester, M13~9PL, UK}
\email{garkusha@imi.ras.ru}
\urladdr{homotopy.nm.ru}

\author{Mike Prest}
\address{School of Mathematics, University of Manchester, Oxford Road, Manchester, M13~9PL, UK}
\date{October 23, 2006}
\email{mprest@manchester.ac.uk}


\urladdr{www.maths.man.ac.uk/$\sim$mprest}


\begin{abstract}
Given a commutative ring $R$ (respectively a positively graded
commutative ring $A=\ps_{j\geq 0}A_j$ which is finitely generated as
an $A_0$-algebra), a bijection between the torsion classes of finite
type in $\Rfp$ (respectively tensor torsion classes of finite type
in $\QGr A$) and the set of all subsets $Y\subseteq\spec R$
(respectively $Y\subseteq\Proj A$) of the form
$Y=\bigcup_{i\in\Omega}Y_i$, with $\spec R\setminus Y_i$
(respectively $\Proj A\setminus Y_i$) quasi-compact and open for all
$i\in\Omega$, is established. Using these bijections, there are
constructed isomorphisms of ringed spaces
   $$(\spec R,\cc O_{R})\lra{\sim}(\spec(\Rfp),\cc O_{\Rfp})$$
and
   $$(\Proj A,\cc O_{\Proj A})\lra{\sim}(\spec(\QGr A),\cc O_{\QGr A}),$$
where $(\spec(\Rfp),\cc O_{\Rfp})$ and $(\spec(\QGr A),\cc O_{\QGr
A})$ are ringed spaces associated to the lattices $L_{\serre}(\Rfp)$
and $L_{\serre}(\QGr A)$ of torsion classes of finite type. Also, a
bijective correspondence between the thick subcategories of perfect
complexes $\perf(R)$ and the torsion classes of finite type in
$\Rfp$ is established.
\end{abstract}

\maketitle
\tableofcontents

\thispagestyle{empty} \pagestyle{plain}

\section{Introduction}

Non-commutative geometry comes in various flavours.  One is based on
abelian and triangulated categories, the latter being replacements
of classical schemes. This is based on classical results of Gabriel
and later extensions, in particular by Thomason. Precisely,
Gabriel~\cite{Ga} proved that any noetherian scheme $X$ can be
reconstructed uniquely up to isomorphism from the abelian category,
$\Qcoh X$, of quasi-coherent sheaves over $X$. This reconstruction
result has been generalized to quasi-compact schemes by Rosenberg
in~\cite{R}. Based on Thomason's classification theorem,
Balmer~\cite{B1} reconstructs a noetherian scheme $X$ from the
triangulated category of perfect complexes $\perf(X)$. This result
has been generalized to quasi-compact, quasi-separated schemes by
Buan-Krause-Solberg~\cite{BKS}.

In this paper we reconstruct affine and projective schemes from
appropriate abelian categories. Our approach, similar to that used
in~\cite{GP,GP1}, is different from Rosenberg's~\cite{R} and
less abstract. Moreover, some results of the paper are of
independent interest.

Let $\Rfp$ (respectively $\QGr A$) denote the category of
$R$-modules (respectively graded $A$-modules modulo torsion modules)
with $R$ (respectively $A=\ps_{n\geq 0}A_n$) a commutative ring
(respectively a commutative graded ring). We first demonstrate the
following result (cf.~\cite{GP,GP1}).

\begin{theo}[Classification]
Let $R$ (respectively $A$) be a commutative ring (respectively
commutative graded ring which is finitely generated as an
$A_0$-algebra). Then the maps
   $$V\mapsto\cc S=\{M\in\Rfp\mid\supp_R(M)\subseteq V\},\quad\cc S\mapsto V=\bigcup_{M\in\cc S}\supp_R(M)$$
and
   $$V\mapsto\cc S=\{M\in\QGr A\mid\supp_A(M)\subseteq V\},\quad\cc S\mapsto V=\bigcup_{M\in\cc S}\supp_A(M)$$
induce bijections between
\begin{enumerate}
 \item the set of all subsets $V\subseteq\spec R$ (respectively $V\subseteq\Proj A$) of the form
       $V=\bigcup_{i\in \Omega} Y_i$ with $\spec R\setminus Y_i$ (respectively $\Proj A\setminus Y_i$)
       quasi-compact and open for all $i\in \Omega$,
 \item the set of all torsion classes of finite type in $\Rfp$ (respectively tensor torsion
classes of finite type in $\QGr A$).
\end{enumerate}
\end{theo}

This theorem says that $\spec R$ and $\Proj A$ contain all the
information about finite localizations in $\Rfp$ and $\QGr A$
respectively. The next result says that there is a 1-1
correspondence between the finite localizations in $\Rfp$ and the
triangulated localizations in $\perf(R)$ (cf.~\cite{Ho,GP}).

\begin{theo}
Let $R$ be a commutative ring. The map
   $$\cc S\mapsto\mathcal T=\{X\in\perf(R)\mid H_n(X)\in\cc S\textrm{ for all }
       n\in\mathbf Z\}$$
induces a bijection between
\begin{enumerate}

\item the set of all torsion classes of finite type in $\Rfp$,

\item the set of all thick subcategories of $\perf(R)$.
\end{enumerate}
\end{theo}

Following Buan-Krause-Solberg~\cite{BKS} we consider the lattices
$L_{\serre}(\Rfp)$ and $L_{\serre}(\QGr A)$ of (tensor) torsion
classes of finite type in $\Rfp$ and $\QGr A$, as well as their
prime ideal spectra $\spec(\Rfp)$ and $\spec(\QGr A)$. These spaces
come naturally equipped with sheaves of rings $\cc O_{\Rfp}$ and
$\cc O_{\QGr A}$. The following result says that the schemes $(\spec
R,\cc O_R)$ and $(\Proj A,\cc O_{\Proj A})$ are isomorphic to
$(\spec(\Rfp),\cc O_{\Rfp})$ and $(\spec(\QGr A),\cc O_{\QGr A})$
respectively.

\begin{theo}[Reconstruction]
Let $R$ (respectively $A$) be a commutative ring (respectively
commutative graded ring which is finitely generated as an
$A_0$-algebra). Then there are natural isomorphisms of ringed spaces
   $$(\spec R,\cc O_{R})\lra{\sim}(\spec(\Rfp),\cc O_{\Rfp})$$
and
   $$(\Proj A,\cc O_{\Proj A})\lra{\sim}(\spec(\QGr A),\cc O_{\QGr A}).$$
\end{theo}

\vskip12pt\noindent \textbf{Acknowledgement.} This paper was written
during the visit of the first author to the University of Manchester
supported by the MODNET Research Training Network in Model Theory.
He would like to thank the University for the kind hospitality.

\section{Torsion classes of finite type}

We refer the reader to the Appendix for necessary facts about
localization and torsion classes in Grothendieck categories.

\begin{prop}\label{aaa}
Assume that $\ff B$ is a set of finitely generated ideals of a
commutative ring $R$. The set of those ideals which contain a finite
products of ideals belonging to $\ff B$ is a Gabriel filter of
finite type.
\end{prop}

\begin{proof}
See~\cite[VI.6.10]{St}.
\end{proof}

Given a module $M$, we denote by $\supp_R(M)=\{P\in\spec R\mid
M_P\neq 0\}$.  Here $M_P$ denotes the localization of $M$ at $P$,
that is, the module of fractions $M[(R\setminus P)^{-1}]$. Note that
$V(I)=\{ P\in \spec R \mid I\leq P\}$ is equal to $\supp_R(R/I)$
for every ideal $I$ and
   $$\supp_R(M)=\bigcup_{x\in M}V(\ann_R(x)),\quad M\in\Rfp.$$

Recall from~\cite{Hoc} that a topological space is {\it spectral\/}
if it is $T_0$, quasi-compact, if the quasi-compact open subsets are
closed under finite intersections and form an open basis, and if
every non-empty irreducible closed subset has a generic point. Given
a spectral topological space, $X$, Hochster \cite{Hoc} endows the
underlying set with a new, ``dual", topology, denoted $X^*$, by
taking as open sets those of the form $Y=\bigcup_{i\in\Omega}Y_i$
where $Y_i$ has quasi-compact open complement $X\setminus Y_i$ for
all $i\in\Omega$. Then $X^*$ is spectral and $(X^*)^*=X$ (see
\cite[Prop.~8]{Hoc}). The spaces, $X$, which we shall consider are
not in general spectral; nevertheless we make the same definition
and denote the space so obtained by $X^*$.

Given a commutative ring $R$, every closed subset of $\spec R$ with
quasi-compact complement has the form $V(I)$ for some finitely
generated ideal, $I$, of $R$ (see~\cite[Chpt.~1, Ex.~17(vii)]{AM}).
Therefore a subset of $\spec^*R$ is open \ifff it is of the form
$\bigcup_\lambda V(I_\lambda)$ with each $I_\lambda$ finitely
generated. Notice that $V(I)$ with $I$ a non-finitely generated
ideal is not open in $\spec^*R$ in general. For instance
(see~\cite[3.16.2]{T}),  let $R=\bb C[x_1,x_2,\ldots]$ and
$m=(x_1,x_2,\ldots)$. It is clear that $V(m)=\{m\}$ is not open in
$\spec^*\bb C[x_1,x_2,\ldots]$.

For definitions of terms used in the next result see the Appendix to this paper.

\begin{thm}[Classification]\label{tol}
Let $R$ be a commutative ring. There are bijections between
\begin{enumerate}
 \item the set of all open subsets $V\subseteq\spec^*R$,
 \item the set of all Gabriel filters $\ff F$ of finite type,
 \item the set of all torsion classes $\cc S$ of finite type in $\Rfp$.
\end{enumerate}
These bijections are defined as follows:
   \begin{align*}
   V\mapsto &
     \left\{
      \begin{array}{rcl}
       \ff F_V&=&\{I\subset R\mid V(I)\subseteq V\}\\
      \cc S_V&=&\{M\in\Rfp\mid\supp_R(M)\subseteq V\}
      \end{array}
    \right.\\
   \ff F\mapsto &
     \left\{
      \begin{array}{rcl}
       V_{\ff F}&=&\bigcup_{I\in\ff F}V(I)\\
       \cc S_{\ff F}&=&\{M\in\Rfp\mid\ann_R(x)\in\ff F\textrm{ for every $x\in M$}\}
      \end{array}
    \right.\\
   \cc S\mapsto &
    \left\{
      \begin{array}{rcl}
      \ff F_{\cc S}&=&\{I\subset R\mid R/I\in\cc S\}\\
      V_{\cc S}&=&\bigcup_{M\in\cc S}\supp_R(M)
      \end{array}
    \right.
   \end{align*}
\end{thm}

\begin{proof}
The bijection between Gabriel filters of finite type and torsion
classes of finite type is a consequence of a theorem of Gabriel
(see, e.g.,~\cite[5.8]{G}).

Let $\ff F$ be a Gabriel filter of finite type.
 Then the set $\Lambda_{\ff F}$ of finitely generated ideals
$I$ belonging to $\ff F$ is a filter basis for $\ff F$.
Therefore
   $V_{\ff F}=\bigcup_{I\in\Lambda_{\ff F}}V(I)$
is open in $\spec^*R$.

Now let $V$ be an open subset of $\spec^*R$. Let $\Lambda$ denote
the set of finitely generated ideals $I$ such that $V(I)\subseteq
V$. By definition of the topology $V=\bigcup_{I\in\Lambda}V(I)$ and
$I_1\cdots I_n\in\Lambda$ for any $I_1,\ldots, I_n\in\Lambda$. We
denote by $\ff F'_V$ the set of ideals $I\subset R$ such that
$I\supseteq J$ for some $J\in\Lambda$. By Proposition~\ref{aaa} $\ff
F'_V$ is a Gabriel filter of finite type. Clearly, $\ff
F'_V\subset\ff F_V=\{I\subset R\mid V(I)\subseteq V\}$. Suppose
$I\in\ff F_V\setminus\ff F'_V$; by~\cite[VI.6.13-15]{St} (cf. the
proof of Theorem~\ref{grad}) there exists a prime ideal $P\in V(I)$
such that $P\not\in\ff F_V'$. But $V(I)\subseteq V$ and therefore
$P\supset J$ for some $J\in\Lambda$, so $P\in\ff F_V'$, a
contradiction. Thus $\ff F_V'=\ff F_V$.

Clearly, $V=V_{\ff F_V}$ for every open subset $V\subseteq\spec^*R$.
Let $\ff F$ be a Gabriel filter of finite type and $I\in\ff F$.
Clearly $\ff F\subset\ff F_{V_{\ff F}}$ and, as above, there is no
ideal belonging to $\ff F_{V_{\ff F}}\setminus\ff F$. We have shown
the bijection between the sets of all Gabriel filters of finite type
and all open subsets in $\spec^*R$. The description of the bijection
between the set of torsion classes of finite type and the set of
open subsets in $\spec^*R$ is now easily checked.
\end{proof}

\section{The fg-topology}

Let $\inj R$ denote the set of isomorphism classes of indecomposable
injective modules. Given a finitely generated ideal $I$ of $R$, we
denote by $\cc S_I$ the torsion class of finite type corresponding
to the Gabriel filter of finite type having $\{I^n\}_{n\geq 1}$ as a
basis (see Proposition~\ref{aaa} and Theorem~\ref{tol}). Note that a
module $M$ has $\cc S_I$-torsion \ifff every element $x\in M$ is
annihilated by some power $I^{n(x)}$ of the ideal $I$. Let us set
   $$D^{\rm fg}(I):=\{E\in\inj R\mid E \textrm{ is $\cc S_I$-torsion
   free}\},\quad V^{\rm fg}(I):=\inj R\setminus D^{\rm fg}(I)$$
(``fg'' referring to this topology being defined using only finitely
generated ideals).

Let $E$ be any indecomposable injective $R$-module. Set $P=P(E)$ to
be the sum of annihilator ideals of non-zero elements, equivalently
non-zero submodules, of $E$.  Since $E$ is uniform the set of
annihilator ideals of non-zero elements of $E$ is closed under
finite sum. It is easy to check (\cite[9.2]{Pr2}) that $P(E)$ is a
prime ideal and $P(E_P)=P$. Here $E_P$ stands for the injective hull
of $R/P$. There is an embedding
   $$\alpha: \spec R\to\inj R,\quad P\mapsto E_P,$$
which need not be surjective. We shall identify $\spec R$
with its image in $\inj R$.

If $P$ is a prime ideal of a commutative ring $R$ its complement in
$R$ is a multiplicatively closed set $S$. Given a module $M$ we
denote the module of fractions $M[S^{-1}]$ by $M_P$. There is a
corresponding Gabriel filter
   $$\ff F^P=\{I\mid P\notin V(I)\}.$$
Clearly, $\ff F^P$ is of finite type. The $\ff F^P$-torsion modules
are characterized by the property that $M_P=0$
(see~\cite[p.~151]{St}).

More generally, let $\cc P$ be a subset of $\spec R$. To $\cc P$ we
associate a Gabriel filter
   $$\ff F^{\cc P}=\bigcap_{P\in\cc P}\ff F^P=\{I\mid\cc P\cap V(I)=\emptyset\}.$$
The corresponding torsion class consists of all modules $M$ with
$M_P=0$ for all $P\in\cc P$.

Given a family of injective $R$-modules $\cc E$, denote by $\ff
F_{\cc E}$ the Gabriel filter determined by $\cc E$. By definition,
this corresponds to the localizing subcategory $\cc S_{\cc
E}=\{M\in\Rfp\mid\Hom_R(M,E)=0 \textrm{ for all $E\in\cc E$}\}$.

\begin{prop}\label{iii}
A Gabriel filter $\ff F$ is of finite type \ifff it is of the form
$\ff F^{\cc P}$ with ${\cc P}$ a closed set in $\spec^*R$. Moreover,
$\ff F^{\cc P}$ is determined by $\cc E_{\cc P}=\{E_P\mid P\in\cc
P\}$: $\ff F^{\cc P} =\{ I\mid \Hom_R(R/I,\cc E_{\cc P})=0\}$.
\end{prop}

\begin{proof}
This is a consequence of Theorem~\ref{tol}.
\end{proof}

\begin{prop}\label{jjj}
Let $\cc P$ be the closure of $P$ in $\spec^*R$. Then $\cc
P=\{Q\in\spec R\mid Q\subseteq P\}$.  Also $\ff F^{P}=\ff F^{\cc P}$.
\end{prop}

\begin{proof}
This is direct from the definition of the topology.
\end{proof}

Recall that for any ideal $I$ of a ring, $R$, and $r\in R$ we have
an isomorphism $R/(I:r)\cong(rR+I)/I$, where $(I:r)=\{s\in R\mid
rs\in I\}$, induced by sending $1+(I:r)$ to $r+I$.

\begin{prop}\label{ppp}
Let $E$ be an indecomposable injective module and let $P(E)$ be the
prime ideal defined before. Let $I$ be a finitely generated ideal of
$R$.  Then $E\in V^{\rm fg}(I)$ \ifff $E_{P(E)}\in V^{\rm fg}(I)$.
\end{prop}

\begin{proof} Let $I$ be such that $E=E(R/I)$. For each $r\in R\setminus
I$ we have, by the remark just above, that the annihilator of
$r+I\in E $ is $(I:r)$ and so, by definition of $P(E), $ we have
$(I:r)\leq P(E)$. The natural projection $(rR+I)/I\cong
R/(I:r)\longrightarrow R/P(E)$ extends to a morphism from $E$ to
$E_{P(E)}$ which is non-zero on $r+I$. Forming the product of these
morphisms as $r$ varies over $R\setminus I$, we obtain a morphism
from $E$ to a product of copies of $E_{P(E)}$ which is monic on
$R/I$ and hence is monic. Therefore $E$ is a direct summand of a
product of copies of $E_{P(E)}$ and so $E\in V^{\rm fg}(J)$ implies
$E_{P(E)}\in V^{\rm fg}(J)$, where $J$ is a finitely generated
ideal.

Now, $E_{P(E)}\in V^{\rm fg}(I)$, where $I$ is a finitely generated
ideal, means that there is a non-zero morphism
$f:R/I^n\longrightarrow E_{P(E)}$ for some $n$. Since $R/P(E)$ is
essential in $E_{P(E)}$ the image of $f$ has non-zero intersection
with $R/P(E)$ so there is an ideal $J$, without loss of generality
finitely generated, with $I^n<J\leq R$, $J/I^n$ a cyclic module, and
such that the restriction, $f'$, of $f$ to $J/I^n$ is non-zero (and
the image is contained in $R/P(E)$). Since $J/I^n$ is a cyclic $\cc
S_I$-torsion module, there is an epimorphism $g:R/I^m\to J/I^n$ for
some $m$. By construction, $R/P(E)=\varinjlim R/I_\lambda$, where
$I_\lambda$ ranges over the annihilators of non-zero elements of
$E$. Since $R/I^m$ is finitely presented, $0\ne f'g$ factorises
through one of the maps $R/I_\lambda\longrightarrow R/P(E)$.  In
particular, there is a non-zero morphism $R/I^m\longrightarrow E$
showing that $E\in V^{\rm fg}(I)$, as required.
\end{proof}

Given a module $M$, we set
   $$[M]:=\{E\in\inj R\mid\Hom_R(M,E)=0\},\quad (M):=\inj R\setminus[M].$$

\begin{rem}\label{ooo}
For any finitely generated ideal $I$ we have: $D^{\rm
fg}(I)\cap\spec R=D(I)$ and $V^{\rm fg}(I)\cap\spec R=V(I)$.
Moreover, $D^{\rm fg}(I)=[R/I]$.
\end{rem}

If $I,J$ are finitely generated ideals, then $D(IJ)=D(I)\cap D(J)$.
It follows from Proposition~\ref{ppp} and Remark~\ref{ooo} that
$D^{\rm fg}(I)\cap D^{\rm fg}(J)=D^{\rm fg}(IJ)$. Thus the sets
$D^{\rm fg}(I)$ with $I$ running over finitely generated ideals form
a basis for a topology on $\inj R$ which we call the {\it fg-ideals
topology}. This topological space will be denoted by $\injfg R$.
Observe that if $R$ is coherent then the fg-topology equals the
Zariski topology on $\inj R$ (see~\cite{Pr2,GP}). The latter
topological space is defined by taking the $[M]$ with $M$ finitely
presented as a basis of open sets.

\begin{thm} \textrm{\em(cf.~Prest~\cite[9.6]{Pr2})} \label{ddd}
Let $R$ be a commutative ring, let $E$ be an indecomposable
injective module and let $P(E)$ be the prime ideal defined before.
Then $E$ and $E_{P(E)}$ are topologically indistinguishable in
$\injfg R$.
\end{thm}

\begin{proof}
This follows from Proposition~\ref{ppp} and Remark~\ref{ooo}.
\end{proof}

\begin{thm} \textrm{\em(cf.~Garkusha-Prest~\cite[Thm.~A]{GP})} \label{bbb}
Let $R$ be a commutative ring. The space $\spec R$ is dense and a
retract in $\injfg R$. A left inverse to the embedding $\spec
R\hookrightarrow\injfg R$ takes an indecomposable injective module
$E$ to the prime ideal $P(E)$. Moreover, $\injfg R$ is
quasi-compact, the basic open subsets $D^{\rm fg}(I)$, with $I$
finitely generated, are quasi-compact, the intersection of two
quasi-compact open subsets is quasi-compact, and every non-empty
irreducible closed subset has a generic point.
\end{thm}

\begin{proof}
For any finitely generated ideal $I$ we have
   $$D^{\rm fg}(I)\cap\spec R=D(I)$$
(see Remark~\ref{ooo}). From this relation and Theorem~\ref{ddd}
it follows that $\spec R$ is dense in $\injfg R$ and that
$\alpha:\spec R\to\injfg R$ is a continuous map.

One may check (see \cite[9.2]{Pr2}) that
   $$\beta:\injfg R\to\spec R,\ \ \ E\mapsto P(E),$$
is left inverse to $\alpha$. Remark~\ref{ooo} implies that
$\beta$ is continuous. Thus $\spec R$ is a retract of
$\injfg R$.

Let us show that each basic open set $D^{\rm fg}(I)$ is
quasi-compact (in particular $\injfg R=D^{\rm fg}(R)$ is quasi-compact). Let
$D^{\rm fg}(I)=\bigcup_{i\in\Omega}D^{\rm fg}(I_i)$ with each $I_i$
finitely generated. It follows from Remark~\ref{ooo} that
$D(I)=\bigcup_{i\in\Omega}D(I_i)$. Since $I$ is finitely generated,
$D(I)$ is quasi-compact in $\spec R$ by~\cite[Chpt.~1,
Ex.~17(vii)]{AM}. We see that $D(I)=\bigcup_{i\in\Omega_0}D(I_i)$
for some finite subset $\Omega_0\subset\Omega$.

Assume $E\in D^{\rm fg}(I)\setminus\bigcup_{i\in\Omega_0}D^{\rm
fg}(I_i)$. It follows from Theorem~\ref{ddd} that $E_{P(E)}\in
D^{\rm fg}(I)\setminus\bigcup_{i\in\Omega_0}D^{\rm fg}(I_i)$. But
$E_{P(E)}\in D^{\rm fg}(I)\cap\spec
R=D(I)=\bigcup_{i\in\Omega_0}D(I_i)$, and hence it is in
$D(I_{i_0})=D^{\rm fg}(I_{i_0})\cap\spec R$ for some
$i_0\in\Omega_0$, a contradiction. So $D^{\rm fg}(I)$ is
quasi-compact. It also follows that the intersection $D^{\rm
fg}(I)\cap D^{\rm fg}(J)=D^{\rm fg}(IJ)$ of two quasi-compact open
subsets is quasi-compact.  Furthermore, every quasi-compact open
subset in $\injfg R$ must therefore have the form $D^{\rm fg}(I)$
with $I$ finitely generated.

Finally, it follows from Remark~\ref{ooo} and Theorem~\ref{ddd}
that a subset $V$ of $\injfg R$ is Zariski-closed and irreducible
\ifff there is a prime ideal $Q$ of $R$ such that $V=\{E\mid
P(E)\geq Q\}$. This obviously implies that the point $E_Q\in V$ is
generic.
\end{proof}

Notice that $\injfg R$ is not a spectral space in general, for it
is not necessarily $T_0$.

\begin{lem}\label{lll}
Let the ring $R$ be commutative. Then the maps
   $$\spec^* R\supseteq V\bl\phi\mapsto\cc Q_{V}=\{E\in\inj R\mid P(E)\in V\}$$
and
   $$(\injfg R)^*\supseteq\cc Q\bl\psi\mapsto V_{\cc Q}=\{P(E)\in\spec^* R\mid E\in\cc Q\}=
     \cc Q\cap\spec^*R$$
induce a 1-1 correspondence between the lattices of open sets of
$\spec^* R$ and those of $(\injfg R)^*$.
\end{lem}

\begin{proof}
First note that $E_P\in\cc Q_{V}$ for any $P\in V$
(see~\cite[9.2]{Pr2}). Let us check that $\cc Q_{V}$ is an open set
in $(\injfg R)^*$. Every closed subset of $\spec R$ with
quasi-compact complement has the form $V(I)$ for some finitely
generated ideal, $I$, of $R$ (see~\cite[Chpt.~1, Ex.~17(vii)]{AM}),
so there are finitely generated ideals $I_\lambda\subseteq R$ such
that $V=\bigcup_\lambda V(I_\lambda)$. Since the points $E$ and
$E_{P(E)}$ are,  by Theorem~\ref{ddd}, indistinguishable in $(\injfg
R)^*$ we see that $\cc Q_{V}=\bigcup_\lambda V^{\rm fg}(I_\lambda)$,
hence this set is open in $(\injfg R)^*$.

The same arguments imply that $V_{\cc Q}$ is open in $\spec^*R$. It
is now easy to see that $V_{\cc Q_{V}}=V$ and $\cc Q_{V_{\cc Q}}=\cc
Q$.
\end{proof}

\section{Torsion classes and thick subcategories}

We shall write $L(\spec^*R)$, $L((\injfg R)^*)$,
$L_{\thick}(\perf(R))$, $L_{\serre}(\Rfp)$ to denote:
\begin{itemize}
\item{$\diamond$} the lattice of all open subsets of $\spec^*R$,
\item{$\diamond$} the lattice of all open subsets of $(\injfg
R)^*$,

\item{$\diamond$} the lattice of all thick subcategories of
$\perf(R)$,

\item{$\diamond$} the lattice of all torsion classes of finite type in $\Rfp$,
ordered by inclusion.
\end{itemize}
(A thick subcategory is a triangulated subcategory closed under
direct summands).

Given a perfect complex $X\in\perf(R)$ denote by
$\supp(X)=\{P\in\spec R\mid X\otimes^L_RR_P\ne 0\}$. It is easy to
see that
   $$\supp(X)=\bigcup_{n\in\bb Z}\supp_R(H_n(X)),$$
where $H_n(X)$ is the $n$th homology group of $X$.

\begin{thm}[Thomason~\cite{T}]\label{nn}
Let $R$ be a commutative ring. The assignments
   $$\cc T\in L_{\thick}(\perf(R))\bl\mu\mapsto\bigcup_{X\in\cc T}\supp(X)$$
and
   $$V\in L(\spec^*R)\bl\nu\mapsto
     \{X\in\perf(R)\mid \supp(X)\subseteq V\}$$
are mutually inverse lattice isomorphisms.
\end{thm}

Given a subcategory $\cc X$ in $\Rfp$, we may consider the smallest
torsion class of finite type in $\Rfp$ containing $\cc X$. This
torsion class we denote by
   $$\surd\cc X=\bigcap\{\cc S\subseteq\Rfp\mid\cc S\supseteq\cc X\textrm{ is a
   torsion class of finite type}\}.$$

\begin{thm}\textrm{\em(cf. Garkusha-Prest~\cite[Thm.~C]{GP})} \label{qqq}
Let $R$ be a commutative ring. There are bijections between
\begin{itemize}

\item{$\diamond$} the set of all open subsets $Y\subseteq(\injfg
R)^*$,
\item{$\diamond$} the set of all torsion classes of finite type in $\Rfp$,

\item{$\diamond$} the set of all thick subcategories of $\perf(R)$.
\end{itemize}
These bijections are defined as follows:
   \begin{align*}
   Y\mapsto &
     \left\{
      \begin{array}{rcl}
       \cc S&=&\{M\mid(M)\subseteq Y\}\\
       \cc T&=&\{X\in\perf(R)\mid (H_n(X))\subseteq Y\textrm{ for all }
       n\in\bb Z\}
      \end{array}
    \right.\\
   \cc S\mapsto &
     \left\{
      \begin{array}{rcl}
       Y&=&\bigcup_{M\in\mathcal S}(M)\\
       \cc T&=&\{X\in\perf(R)\mid H_n(X)\in\cc S\textrm{ for all }
       n\in\bb Z\}
      \end{array}
    \right.\\
   \cc T\mapsto &
     \left\{
      \begin{array}{rcl}
       Y&=&\bigcup_{X\in\mathcal T,n\in\bb Z}(H_n(X))\\
       \cc S&=&\surd\{H_n(X)\mid X\in\cc T, n\in\bb Z\}
      \end{array}
    \right.
   \end{align*}
\end{thm}

\begin{proof}
That $\cc S_Y=\{M\mid(M)\subseteq Y\}$ is a torsion class follows
because it is defined as the class of modules having no non-zero
morphism to a family of injective modules, $\cc E:=\inj R\setminus
Y$. By Lemma~\ref{lll}, $\cc E\cap\spec^*R=U$ is a closed set in
$\spec^*R$, that is $P(E)\in U$ for all $E\in\cc E$. $\cc S_Y$ is
also determined by the family of injective modules $\{E_P\}_{P\in
U}$. Indeed, any $E\in\cc E$ is a direct summand of some power of
$E_{P(E)}$ by the proof of Proposition~\ref{ppp}. Therefore
$\Hom_R(M,E_{P(E)})=0$ implies $\Hom_R(M,E)=0$. By
Proposition~\ref{iii} $\cc S_Y$ is of finite type. Conversely, given
a torsion class of finite type $\cc S$, the set $Y_{\cc
S}=\bigcup_{M\in\mathcal S}(M)$ is plainly open in $(\injfg R)^*$.
Moreover, $\cc S_{Y_{\cc S}}=\cc S$ and $Y=Y_{\cc S_Y}$.

Consider the following diagram:
   $$\xymatrix{L(\spec^*R)\ar@<.6ex>[rr]^\nu\ar@<.6ex>[d]^\phi
               &&L_{\thick}(\perf(R))\ar@<.6ex>[ll]^\mu\ar@<.6ex>[d]^\rho\\
               L((\injfg R)^*)\ar@<.6ex>[rr]^\zeta\ar@<.6ex>[u]^\psi
               &&L_{\serre}(\Rfp),\ar@<.6ex>[ll]^\delta\ar@<.6ex>[u]^\sigma}$$
where $\phi,\psi$ are as in Lemma~\ref{lll}, $\mu,\nu$ are as in
Theorem~\ref{nn} and the remaining maps are the corresponding maps
indicated in the formulation of the theorem. We have $\nu=\mu^{-1}$ by
Theorem~\ref{nn}, $\phi=\psi^{-1}$ by Lemma~\ref{lll}, and
$\zeta=\delta^{-1}$ by the above.

By construction,
   $$\sigma\zeta\phi(V)=\{X\mid\bigcup_{n\in\bb Z}\supp_R(H_n(X))
   \subseteq V\}=\{X\mid\supp(X)\subseteq V\}$$
for all $V\in L(\spec^*R)$. Thus $\sigma\zeta\phi=\nu$. Since
$\zeta,\phi,\nu$ are bijections so is $\sigma$.

On the other hand,
   $$\psi\delta\rho(\cc T)=\bigcup_{X\in\cc T, n\in\bb Z}\supp_R(H_n(X))
   =\bigcup_{X\in\cc T}\supp(X)$$
for any $\cc T\in L_{\thick}(\perf(R))$. We have used here the
relation
   $$\bigcup_{M\in\rho(\cc T)}\supp_R(M)=\bigcup_{X\in\cc T, n\in\bb Z}\supp_R(H_n(X)).$$
One sees that $\psi\delta\rho=\mu$. Since $\delta,\psi,\mu$ are
bijections so is $\rho$. Obviously, $\sigma=\rho^{-1}$ and the
diagram above yields the desired bijective correspondences. The theorem
is proved.
\end{proof}

To conclude this section, we should mention the relation between
torsion classes of finite type in $\Rfp$ and the Ziegler subspace
topology on $\inj R$ (we denote this space by $\injzg R$). The
latter topology arises from Ziegler's work on the model theory of
modules~\cite{Z}. The points of the Ziegler spectrum of $R$ are the
isomorphism classes of indecomposable pure-injective $R$-modules and
the closed subsets correspond to complete theories of modules. It is
well known (see \cite[9.12]{Pr2}) that for every coherent ring $R$
there is a 1-1 correspondence between the open (equivalently closed)
subsets of $\injzg R$ and torsion classes of finite type in $\Rfp$.
However, this is not the case for general commutative rings.

The topology on $\injzg R$ can be defined as follows. Let $\cc M$ be
the set of those modules $M$ which are kernels of homomorphisms
between finitely presented modules; that is $M=\kr(K\lra{f}L)$ with
$K,L$ finitely presented. The sets $(M)$ with $M\in\cc M$ form a
basis of open sets for $\injzg R$. We claim that there is a ring $R$
and a module $M\in\cc M$ such that the intersection
$(M)\cap\spec^*R$ is not open in $\spec^*R$, and hence such that the
open subset $(M)$ cannot correspond to any torsion class of finite
type on $\Rfp$. Such a ring has been pointed out by G.~Puninski.

Let $V$ be a commutative valuation domain with value group
isomorphic to $\Gamma=\ps_{n\in\bb Z}\bb Z$, a $\bb Z$-indexed
direct sum of copies of $\bb Z$. The order on $\Gamma$ is defined as
follows. $(a_n)_{n\in\bb Z}>(b_n)_{n\in\bb Z}$ if $a_i>b_i$ for some
$i$ and $a_k=b_k$ for every $k<i$. Then $J^2=J$ where $J$ is the
Jacobson radical of $V$. Let $r$ be an element with value
$v(r)=(a_n)_{n\in\bb Z}$ where $a_0=1$ and $a_n=0$ for all $n\not=
0$. Consider the ring $R=V/rJ$. Again $J(R)^2=J(R)$. Denoting the
image of $r$ in $R$ by $r'$, note that $\ann_R(r')=J(R)$ which is
not finitely generated, and so $R$ is not coherent by Chase's
Theorem (see~\cite[1.13.3]{St}). Note that $R$ is a local ring and,
as already observed, the simple module  $R/J(R)$ is isomorphic to
$r'R$.  Therefore $R/J(R)=\kr(R\to R/r'R)$. Thus $R/J(R)\in\cc M$.
We have
   $$(R/J(R))\cap\spec^*R=V(J(R))=\{J(R)\}.$$
Suppose $V(J(R))$ is open in $\spec^*R$; then
$V(J(R))=\bigcup_\lambda V(I_\lambda)$ with each $I_\lambda$
finitely generated. Since $J(R)$ is the largest proper ideal each
$V(I_\lambda)$, if non-empty, equals $\{ J(R)\}$.  Therefore
$J(R)=\surd I_\lambda$ for some $\lambda$. But the prime radical of
every finitely generated ideal in $R$ is prime (since $R$ is a
valulation ring) and different from $J(R)$.  To see the latter, we
have, since $I_\lambda$ is finitely generated, that all elements of
$I_\lambda$ have value $>(a'_n)_n$ for some $(a'_n)_n$ with $a'_n=0$
for all $n\leq N$ for some fixed $N$.  (Recall that the valuation
$v$ on $R$ satisfies $v(r+s)\geq {\rm min}\{ v(r), v(s)\}$ and
$v(rs)=v(r)+v(s)$.)  It follows that there is a prime ideal properly
between $I_\lambda$ and $J(R)$.  This gives a contradiction, as
required.

\section{Graded rings and modules}

In this section we recall some basic facts about graded rings and
modules.

\begin{defs}{\rm
A {\it (positively) graded ring\/} is a ring $A$ together with a
direct sum decomposition $A=A_0\ps A_1\ps A_2\ps\cdots$ as abelian
groups, such that $A_iA_j\subset A_{i+j}$ for $i,j\geq 0$. A {\it
homogeneous element\/} of $A$ is simply an element of one of the
groups $A_j$, and a {\it homogeneous ideal\/} of $A$ is an ideal
that is generated by homogeneous elements. A {\it graded
$A$-module\/} is an $A$-module $M$ together with a direct sum
decomposition $M=\ps_{j\in\bb Z}M_j$ as abelian groups, such that
$A_iM_j\subset M_{i+j}$ for $i\geq 0,j\in\bb Z$. One calls $M_j$ the
$j$th {\it homogeneous component of $M$}. The elements $x\in M_j$
are said to be {\it homogeneous (of degree $j$)}.

Note that $A_0$ is a commutative ring with $1\in A_0$, that all
summands $M_j$ are $A_0$-modules, and that $M=\ps_{j\in\bb Z}M_j$ is
a direct sum decomposition of $M$ as an $A_0$-module.

Let $A$ be a graded ring. The {\it category of graded $A$-modules},
denoted by $\Gr A$, has as objects the graded $A$-modules. A {\it
morphism\/} of graded $A$-modules $f:M\to N$ is an $A$-module
homomorphism satisfying $f(M_j)\subset N_j$ for all $j\in\bb Z$. An
$A$-module homomorphism which is a morphism in $\Gr A$ will be
called {\it homogeneous}.

Let $M$ be a graded $A$-module and let $N$ be a submodule of $M$.
Say that $N$ is a {\it graded submodule\/} if it is a graded module
such that the inclusion map is a morphism in $\Gr A$. The graded
submodules of $A$ are called {\it graded ideals}. If $d$ is an
integer the {\it tail\/} $M_{\geq d}$ is the graded submodule of $M$
having the same homogeneous components $(M_{\geq d})_j$ as $M$ in
degrees $j\geq d$ and zero for $j<d$. We also denote the ideal
$A_{\geq 1}$ by $A_+$.
}\end{defs}

For $n\in\bb Z$, $\Gr A$ comes equipped with a shift functor
$M\mapsto M(n)$ where $M(n)$ is defined by $M(n)_j=M_{n+j}$. Then $\Gr A$ is a
Grothendieck category with generating family $\{A(n)\}_{n\in\bb
Z}$. The tensor product for the category of all $A$-modules induces
a tensor product on $\Gr A$: given two graded $A$-modules $M,N$ and
homogeneous elements $x\in M_i,y\in N_j$, set $\deg(x\otimes
y):=i+j$. We define the {\it homomorphism $A$-module\/} $\cc
Hom_A(M,N)$ to be the graded $A$-module which is, in dimension $n\in\bb Z$, the group $\cc
Hom_A(M,N)_n$ of graded $A$-module homomorphisms of
degree $n$, i.e.,
   $$\cc Hom_A(M,N)_n=\Gr A(M,N(n)).$$

We say that a graded $A$-module $M$ is {\it finitely generated\/} if
it is a quotient of a free graded module of finite rank
$\bigoplus_{s=1}^nA(d_s)$ where $d_1,\ldots,d_s\in\bb Z$. Say that $M$ is
{\it finitely presented\/} if there is an exact sequence
   $$\bigoplus_{t=1}^mA(c_t)\to\bigoplus
   _{s=1}^nA(d_s)\to M\to 0.$$
The full subcategory of finitely presented graded modules will be
denoted by $\gr A$. Note that any graded $A$-module is a direct
limit of finitely presented graded $A$-modules, and therefore $\Gr
A$ is a locally finitely presented Grothendieck category.

Let $E$ be any indecomposable injective graded $A$-module (we remind
the reader that the corresponding ungraded module, $\bigoplus_nE_n$,
need not be injective in the category of ungraded $A$-modules). Set
$P= P(E)$ to be the sum of the annihilator ideals $\ann_A(x)$ of
non-zero homogeneous elements $x\in E$. Observe that each ideal
$\ann_A(x)$ is homogeneous. Since $E$ is uniform the set of
annihilator ideals of non-zero homogeneous elements of $E$ is upwards closed
so the only issue is whether the sum, $P(E)$, of
them all is itself one of these annihilator ideals.

Given a prime homogeneous ideal $P$, we use the notation $E_P$ to
denote the injective hull, $E(A/P)$, of $A/P$. Notice that $E_P$ is
indecomposable. We also denote the set of isomorphism classes of
indecomposable injective graded $A$-modules by $\inj A$.

\begin{lem}\label{9.2}
If $E\in \inj A$ then $P(E)$ is a homogeneous prime ideal. If the
module $E$ has the form $E_P(n)$ for some prime homogeneous ideal
$P$ and integer $n$, then $P=P(E)$.
\end{lem}
\begin{proof}
The proof is similar to that of~\cite[9.2]{Pr2}.
\end{proof}

It follows from the preceding lemma that the map
   $$P\subset A\mapsto E_P\in\inj A$$
from the set of homogeneous prime ideals to $\inj A$ is injective.

A {\it tensor torsion class\/} in $\Gr A$ is a torsion class with
torsion class $\cc S\subset\Gr A$ such that for any $X\in\cc S$ and
any $Y\in\Gr A$ the tensor product $X\otimes Y$ is in $\cc S$.

\begin{lem}\label{ccc1}
Let $A$ be a graded ring. Then a torsion class $\cc S$ is a tensor
torsion class of $\Gr A$ \ifff it is closed under shifts of objects,
i.e.~$X\in\cc S$ implies $X(n)\in\cc S$ for any $n\in\bb Z$.
\end{lem}

\begin{proof}
Suppose that $\cc S$ is a tensor torsion class of $\Gr A$. Then it
is closed under shifts of objects, because $X(n)\cong X\otimes
A(n)$.

Assume the converse. Let $X\in\cc S$ and $Y\in\Gr A$. Then there is
a surjection $\ps_{i\in I}A(i)\bl f\twoheadrightarrow Y$. It follows
that $1_X\otimes f:\ps_{i\in I}X(i)\to X\otimes Y$ is a surjection.
Since each $X(i)$ belongs to $\cc S$ then so does $X\otimes Y$.
\end{proof}

\begin{lem}\label{Gabr}
The map
   $$\cc S\longmapsto\ff F(\cc S)=\{\ff a\subseteq A\mid A/\ff a\in\cc S\}$$
establishes a bijection between the tensor torsion classes in $\Gr
A$ and the sets $\ff F$ of homogeneous ideals satisfying the
following axioms:
\begin{itemize}
\item[$T1.$] $A\in\ff F$;
\item[$T2.$] if $\ff a\in\ff F$ and $a$ is a homogeneous element of
$A$ then $(\ff a:a)=\{x\in A\mid xa\in\ff a\}\in\ff F$;
\item[$T3.$] if $\ff a$ and $\ff b$ are homogeneous ideals of $A$ such that
    $\ff a\in\ff F$ and $(\ff b:a)\in\ff F$
    for every homogeneous element $a\in\ff a$ then $\ff b\in\ff F$.
\end{itemize}
We shall refer to such filters as {\em $t$-filters}. Moreover, $\cc
S$ is of finite type \ifff $\ff F(\cc S)$ has a basis of finitely
generated ideals, that is every ideal in $\ff F(\cc S)$ contains a
finitely generated ideal belonging to $\ff F(\cc S)$. In this case
$\ff F(\cc S)$ will be refered to as a {\em $t$-filter of finite
type}.
\end{lem}

\begin{proof}
It is enough to observe that there is a bijection between the
Gabriel filters on the family $\{A(n)\}_{n\in\bb Z}$ of generators
closed under the shift functor (i.e., if $\ff a$ belongs to the
Gabriel filter then so does $\ff a(n)$ for all $n\in\bb Z$) and the
$t$-filters.
\end{proof}

\begin{prop}\label{eee}
The following statements are true:
\begin{enumerate}
\item Let $\ff F$ be a $t$-filter. If $I,J$ belong to $\ff F$, then
$IJ\in\ff F$.
\item Assume that $\ff B$ is a set of homogeneous finitely generated ideals. The
set $\ff B'$ of finite products of ideals belonging to $\ff B$ is a
basis for a $t$-filter of finite type.
\end{enumerate}
\end{prop}

\begin{proof}
(1). For any homogeneous element $a\in I$ we have $(IJ:a)\supset J$,
so $IJ\in\ff F$ by $T3$ and the fact that every homogeneous ideal
containing an ideal from $\ff F$ must belong to $\ff F$.

(2). We follow~\cite[VI.6.10]{St}. We must check that the set $\ff
F$ of homogeneous ideals containing ideals in $\ff B'$ is a
$t$-filter of finite type. $T1$ is plainly satisfied. Let $a$ be a
homogeneous element in $A$ and $I\in\ff F$. There is an ideal
$I'\in\ff B'$ contained in $I$. Then $(I:a)\supset I'$ and therefore
$(I:a)\in \ff F$, hence $T2$ is satisfied as well.

Next we verify that $\ff F$ satisfies $T3$. Suppose that $I$ is a
homogeneous ideal and there exists $J\in\ff F$ such that
$(I:a)\in\ff F$ for every homogeneous element $a\in J$. We may
assume that $J\in\ff B'$. Let $a_1,\ldots,a_n$ be generators of
$J$. Then $(I:a_i)\in\ff F$, $i\leq n$, and $(I:a_i)\supset J_i$ for
some $J_i\in\ff B'$. It follows that $a_iJ_i\subset I$ for each $i$,
and hence $JJ_1\cdots J_n\subset J(J_1\cap\ldots\cap J_n)\subset I$,
so $I\in\ff F$.
\end{proof}

\section{Torsion modules and the category $\QGr A$}

In this section we introduce the category $\QGr A$, which is
analogous to the category of quasi-coherent sheaves on a projective
variety. The non-commutative analog of the category $\QGr A$ plays a
prominent role in ``non-commutative projective geometry" (see, e.g.,
\cite{AZ,S,V}).

Recall that the projective scheme $\Proj A$ is a topological space
whose points are the homogeneous prime ideals not containing $A_+$.
The topology of $\Proj A$ is defined by taking the closed sets to be
the sets of the form $V(I)=\{P\in\Proj A\mid P\supseteq I\}$ for
$I$ a homogeneous ideal of $A$. We set $D(I):=\Proj A\setminus
V(I)$. The space $\Proj A$ is spectral and the quasi-compact open
sets are those of the form $D(I)$ with $I$ finitely generated
(see, e.g.,~\cite[5.1]{GP1}).

{\bf In the remainder of this section the homogeneous ideal
$A_+\subset A$ is assumed to be finitely generated.} This is
equivalent to assuming that $A$ is a finitely generated
$A_0$-algebra. Let $\Tors A$ denote the tensor torsion class of
finite type corresponding to the family of homogeneous finitely
generated ideals $\{A^n_+\}_{n\geq 1}$ (see Proposition~\ref{eee}).
We refer to the objects of $\Tors A$ as {\it torsion graded
modules}.

Let $\QGr A=\Gr A/\Tors A$. Let $Q$ denote the quotient functor $\Gr
A\to\QGr A$. We shall identify $\QGr A$ with the full subcategory of
$\Tors$-closed modules. The shift functor $M\mapsto M(n)$ defines a
shift functor on $\QGr A$ for which we shall use the same notation.
Observe that $Q$ commutes with the shift functor. Finally we shall
write $\cc O=Q(A)$. Note that $\QGr A$ is a locally finitely
generated Grothendieck category with the family, $\{Q(M)\}_{M\in\gr
A}$, of finitely generated generators (see~\cite[5.8]{G}).

The tensor product in $\Gr A$ induces a tensor product in $\QGr A$,
denoted by $\boxtimes$. More precisely, one sets
   $$X\boxtimes Y:=Q(X\otimes Y)$$
for any $X,Y\in\QGr A$.

\begin{lem}\label{d2}
Given $X,Y\in\Gr A$ there is a natural isomorphism in $\QGr A$:
$Q(X)\boxtimes Q(Y)\cong Q(X\otimes Y)$. Moreover, the functor
$-\boxtimes Y:\QGr A\to\QGr A$ is right exact and preserves direct
limits.
\end{lem}

\begin{proof}
See \cite[4.2]{GP1}.
\end{proof}

As a consequence of this lemma we get an isomorphism $X(d)\cong\cc
O(d)\boxtimes X$ for any $X\in\QGr A$ and $d\in\bb Z$.

The notion of a tensor torsion class of $\QGr A$ (with respect to
the tensor product $\boxtimes$) is defined analogously to that in
$\Gr A$. The proof of the next lemma is like that of
Lemma~\ref{ccc1} (also use Lemma~\ref{d2}).

\begin{lem}\label{ccc2}
A torsion class $\cc S$ is a tensor torsion class of $\QGr A$ \ifff it is closed
under shifts of objects, i.e.~$X\in\cc S$ implies $X(n)\in\cc S$ for
any $n\in\bb Z$.
\end{lem}

Given a prime ideal $P\in\Proj A$ and a graded module $M$, denote by
$M_P$ the homogeneous localization of $M$ at $P$. If $f$ is a
homogeneous element of $A$, by $M_f$ we denote the localization of
$M$ at the multiplicative set $S_f=\{f^n\}_{n\geq 0}$.

\begin{lem}\label{g111}
If $T$ is a torsion module then $T_P=0$ and $T_f=0$ for any
$P\in\Proj A$ and $f\in A_+$. As a consequence, $M_P\cong Q(M)_P$
and $M_f\cong Q(M)_f$ for any $M\in\Gr A$.
\end{lem}
\begin{proof}
See \cite[5.5]{GP1}.
\end{proof}

Denote by $L_{\serre}(\Gr A,\Tors A)$ (respectively $L_{\serre}(\QGr
A)$) the lattice of the tensor torsion classes of finite type in
$\Gr A$ with torsion classes containing $\Tors A$ (respectively the
tensor torsion classes of finite type in $\QGr A$) ordered by
inclusion. The map
   $$\ell:L_{\serre}(\Gr A,\Tors A)\lra{}L_{\serre}(\QGr A),\quad\cc S\longmapsto\cc S/\Tors A$$
is a lattice isomorphism, where $\cc S/\Tors A=\{Q(M)\mid M\in\cc
S\}$ (see, e.g.,~\cite[1.7]{G}). We shall consider the map $\ell$ as an
identification.

\begin{thm}[Classification]\label{grad}
Let $A$ be a graded ring which is finitely generated as an
$A_0$-algebra. Then the maps
   $$V\mapsto\cc S=\{M\in\QGr A\mid\supp_A(M)\subseteq V\}\quad\textrm{and}\quad\cc S\mapsto V=\bigcup_{M\in\cc S}\supp_A(M)$$
induce bijections between
\begin{enumerate}
 \item the set of all open subsets $V\subseteq\Proj^*A$,
 \item the set of all tensor torsion classes of finite type in $\QGr A$.
\end{enumerate}
\end{thm}

\begin{proof}
By Lemma~\ref{Gabr} it is enough to show that the maps
   $$V\mapsto\ff F_V=\{I\in A\mid V(I)\subseteq V\}\quad
     \textrm{and}\quad\ff F\mapsto V_{\ff F}=\bigcup_{I\in\ff F}V(I)$$
induce bijections between the set of all open subsets
$V\subseteq\Proj^*A$ and the set of all $t$-filters of finite type
containing $\{A_+^n\}_{n\geq 1}$.

Let $\ff F$ be such a $t$-filter. Then the set $\Lambda_{\ff F}$ of
finitely generated graded ideals $I$ belonging to $\ff F$ is a basis
for $\ff F$. Clearly $V_{\ff F}=\bigcup_{I\in\Lambda_{\ff F}}V(I)$,
so $V_{\ff F}$ is open in $\Proj^*A$.

Now let $V$ be an open subset of $\Proj^*A$. Let $\Lambda$ be the
set of finitely generated homogeneous ideals $I$ such that
$V(I)\subseteq V$. Then $V=\bigcup_{I\in\Lambda}V(I)$ and $I_1\cdots
I_n\in\Lambda$ for any $I_1,\ldots, I_n\in\Lambda$. We denote by
$\ff F'_V$ the set of homogeneous ideals $I\subset A$ such that
$I\supseteq J$ for some $J\in\Lambda$. By Proposition~\ref{eee}(2)
$\ff F'_V$ is a $t$-filter of finite type. Clearly, $\ff
F'_V\subset\ff F_V$. Suppose $I\in\ff F_V\setminus\ff F'_V$.

We can use Zorn's lemma to find an ideal $J\supset I$ which is
maximal with respect to $J\notin\ff F'_V$ (we use the fact that $\ff
F'_V$ has a basis of finitely generated ideals). We claim that $J$
is prime. Indeed, suppose $a,b\in A$ are two homogeneous elements
not belonging to $J$. Then $J+aA$ and $J+bA$ must be members of $\ff
F'_V$, and also $(J+aA)(J+bA)\in\ff F'_V$ by
Proposition~\ref{eee}(1). But $(J+aA)(J+bA)\subset J+abA$, and
therefore $ab\notin J$. We see that $J\in V(I)\subset V$, and hence
$J\in V(I')$ for some $I'\in\Lambda$. But this implies $J\in\ff
F_V'$, a contradiction. Thus $\ff F_V'=\ff F_V$.
Clearly, $V=V_{\ff F_V}$ for every open subset $V\subseteq\Proj^*A$.
Let $\ff F$ be a $t$-filter of finite type and $I\in\ff F$. Then
$I\supset J$ for some $J\in\Lambda_{\ff F}$, and hence $V(I)\subset
V(J)\subset V_{\ff F}$. It follows that $\ff F\subset\ff F_{V_{\ff
F}}$. As above, there is no ideal belonging to $\ff F_{V_{\ff
F}}\setminus\ff F$. We have shown the desired bijection between the
sets of all $t$-filters of finite type and all open subsets in
$\Proj^*A$.
\end{proof}

\section{The prime spectrum of an ideal lattice}

Inspired by recent work of Balmer~\cite{B}, Buan, Krause, and
Solberg~\cite{BKS} introduce the notion of an ideal lattice and
study its prime ideal spectrum. Applications arise from
abelian or triangulated tensor categories.

\begin{defs}[Buan, Krause, Solberg~\cite{BKS}] {\rm
An {\it ideal lattice\/} is by definition a partially ordered set
$L=(L,\leq)$, together with an associative multiplication $L\times
L\to L$, such that the following holds.
\begin{enumerate}
\item[(L1)] The poset $L$ is a {\it complete lattice\/}, that is,
$$\sup A = \bigvee_{a\in A} a\quad\text{and}\quad \inf A = \bigwedge_{a\in A}
a$$ exist in $L$ for every subset $A\subseteq L$.
\item[(L2)] The lattice $L$ is {\it compactly generated\/}, that is,
every element in $L$ is the supremum of a set of compact elements.  (An
element $a\in L$ is {\em compact}, if for all $A\subseteq L$ with
$a\leq \sup A$ there exists some finite $A'\subseteq A$ with
$a\leq\sup A'$.)
\item[(L3)] We have for all $a,b,c\in L$
$$a(b\vee c)=ab\vee ac\quad\text{and}\quad (a\vee b)c=ac\vee bc.$$
\item[(L4)] The element $1=\sup L$ is compact, and $1a=a=a1$ for all $a\in L$.
\item[(L5)] The product of two compact elements is again compact.
\end{enumerate}
A {\it morphism\/} $\phi\colon L\to L'$ of ideal lattices is a map
satisfying
\begin{gather*}\label{eq:mor}
\phi(\bigvee_{a\in A}a)=\bigvee_{a\in A}\phi(a)\quad \text{for}\quad
A\subseteq L, \\
\phi(1)=1\quad\text{and}\quad\phi(ab)=\phi(a)\phi(b)\quad\text{for}\quad
a,b\in L.\notag
\end{gather*}
}\end{defs}

Let $L$ be an ideal lattice. Following~\cite{BKS} we define the
spectrum of prime elements in $L$. An element $p\neq 1$ in $L$ is
 {\em prime} if $ab\leq p$ implies $a\leq p$ or $b\leq p$ for
all $a,b\in L$. We denote by $\spec L$ the set of prime elements in
$L$ and define for each $a\in L$
   $$V(a)=\{p\in\spec L\mid a\leq p\}\quad\text{and}\quad D(a)=\{p\in\spec L\mid a\not\leq p\}.$$
The subsets of $\spec L$ of the form $V(a)$ are closed under forming
arbitrary intersections and finite unions.  More precisely,
   $$V(\bigvee_{i\in\Omega} a_i)=\bigcap_{i\in\Omega} V(a_i)\quad\text{and}\quad
   V(ab)=V(a)\cup V(b).$$
Thus we obtain the {\it Zariski topology\/} on $\spec L$ by
declaring a subset of $\spec L$ to be {\it closed\/} if it is of the
form $V(a)$ for some $a\in L$. The set $\spec L$ endowed with this
topology is called the {\it prime spectrum\/} of $L$.  Note that the
sets of the form $D(a)$ with compact $a\in L$ form a basis of open
sets. The prime spectrum $\spec L$ of an ideal lattice $L$ is
spectral~\cite[2.5]{BKS}.

There is a close relation between spectral spaces and ideal
lattices. Given a topological space $X$, we denote by $L_{\open}(X)$
the lattice of open subsets of $X$ and consider the multiplication
map
   $$L_{\open}(X)\times L_{\open}(
 X)\to L_{\open}(X),\quad (U,V)\mapsto UV=U\cap V.$$
The lattice $L_{\open}(X)$ is complete.

The following result, which appears in~\cite{BKS}, is part of the
Stone Duality Theorem (see, for instance, \cite{Jo}).

\begin{prop}\label{pr:openlattice}
Let $X$ be a spectral space. Then $L_{\open}(X)$ is an ideal
lattice. Moreover, the map
   $$X\to\spec L_{\open}(X),\quad x \mapsto X\setminus \overline{\{x\}},$$
is a homeomorphism.
\end{prop}

We deduce from the classification of torsion classes of finite type
(Theorems~\ref{tol} and \ref{grad}) the following.

\begin{prop}\label{prr}
Let $R$ (respectively $A$) be a ring (respectively graded ring which
is finitely generated as an $A_0$-algebra). Then $L_{\serre}(\Rfp)$
and $L_{\serre}(\QGr A)$ are ideal lattices.
\end{prop}

\begin{proof}
The spaces $\spec R$ and $\Proj A$ are spectral. Thus $\spec^*R$ and
$\Proj^*A$ are spectral, also $L_{\open}(\spec^*R)$ and
$L_{\open}(\Proj^*A)$ are ideal lattices by
Proposition~\ref{pr:openlattice}. By Theorems~\ref{tol} and
\ref{grad} we have isomorphisms $L_{\open}(\spec^*R)\cong
L_{\serre}(\Rfp)$ and $L_{\open}(\Proj^*A)\cong L_{\serre}(\QGr A)$.
Therefore $L_{\serre}(\Rfp)$ and $L_{\serre}(\QGr A)$ are ideal
lattices.
\end{proof}

\begin{cor}\label{prrco}
The points of $\spec L_{\serre}(\Rfp)$ (respectively $\spec
L_{\serre}(\QGr A)$) are the $\cap$-irreducible torsion classes of
finite type in $\Rfp$ (respectively tensor torsion classes of finite
type in $\QGr A$) and the maps
   \begin{gather*}\label{ggg}
    f:\spec^*R\lra{}\spec L_{\serre}(\Rfp),\quad P\longmapsto\cc S_P=\{M\in\Rfp\mid
    M_P=0\}\\
    f:\Proj^*A\lra{}\spec L_{\serre}(\QGr A),\quad P\longmapsto\cc S_P=\{M\in\QGr A\mid
M_P=0\}
   \end{gather*}
are homeomorphisms of spaces.
\end{cor}

\begin{proof}
This is a consequence of Theorems~\ref{tol}, \ref{grad} and
Propositions~\ref{pr:openlattice}, \ref{prr}.
\end{proof}

\section{Reconstructing affine and projective schemes}

We shall write $\spec(\Rfp):=\spec^*L_{\serre}(\Rfp)$ (respectively
$\spec(\QGr A):=\spec^*L_{\serre}(\QGr A)$) and $\supp(M):=\{\cc
P\in\spec(\Rfp)\mid M\not\in\cc P\}$ (respectively $\supp(M):=\{\cc
P\in\spec(\QGr A)\mid M\not\in\cc P\}$) for $M\in\Rfp$ (respectively
$M\in\QGr A$). It follows from Corollary~\ref{prrco} that
   $$\supp_R(M)=f^{-1}(\supp(M))\quad\textrm{(respectively $\supp_A(M)=f^{-1}(\supp(M))$)}.$$

Following \cite{B,BKS}, we define a structure sheaf on $\spec(\Rfp)$
($\spec(\QGr A)$) as follows.  For an open subset $U\subseteq
\spec(\Rfp)$ ($U\subseteq \spec(\QGr A)$), let
   $$\cc S_U=\{M\in\Rfp\ (\QGr A)\mid\supp(M)\cap U=\emptyset\}$$
and observe that $\cc S_U=\{M\mid M_P=0\textrm{ for all $P\in
f^{-1}(U)$}\}$ is a (tensor) torsion class. We obtain a presheaf of
rings on $\spec(\Rfp)$ ($\spec(\QGr A)$) by
   $$U\mapsto\End_{\Rfp/\cc S_U}(R)\quad(\End_{\QGr A/\cc S_U}(\cc O)).$$
If $V\subseteq U$ are open subsets, then the restriction map
   $$\End_{\Rfp/\cc S_U}(R)\to\End_{\Rfp/\cc S_V}(R)\quad(\End_{\QGr A/\cc S_U}(\cc O)
   \to\End_{\QGr A/\cc S_V}(\cc O))$$
is induced by the quotient functor $\Rfp/\cc S_U\to\Rfp/\cc S_V$
($\QGr A/\cc S_U\to\QGr A/\cc S_V$). The sheafification is called
the {\it structure sheaf\/} of $\Rfp$ ($\QGr A$) and is denoted by
$\cc O_{\Rfp}$ ($\cc O_{\QGr A}$). This is a sheaf of commutative
rings by~\cite[XI.2.4]{K}. Next let $\cc P\in\spec(\Rfp)$ and
$P:=f^{-1}(\cc P)$. We have
   $$\cc O_{\Rfp,\cc P}=\lp_{\cc P\in U}\End_{\Rfp/\cc S_U}(R)=\lp_{f\notin P}
   \End_{\Rfp/\cc S_{D(f)}}(R)
     \cong\lp_{f\notin P}R_f=\cc O_{R,P}.$$
Similarly, for $\cc P\in\spec(\QGr A)$ and $P:=f^{-1}(\cc P)$ we
have
   $$\cc O_{\QGr A,\cc P}\cong\cc O_{\Proj A,P}.$$

The next theorem says that the abelian category $\Rfp$ ($\QGr A$)
contains all the necessary information to reconstruct the affine
(projective) scheme $(\spec R,\cc O_{R})$ (respectively $(\Proj
A,\cc O_{\Proj A})$).

\begin{thm}[Reconstruction]\label{coh}
Let $R$ (respectively $A$) be a ring (respectively graded ring which
is finitely generated as an $A_0$-algebra). The maps of
Corollary~\ref{ggg} induce isomorphisms of ringed spaces
   $$f:(\spec R,\cc O_{R})\lra{\sim}(\spec(\Rfp),\cc O_{\Rfp})$$
and
   $$f:(\Proj A,\cc O_{\Proj A})\lra{\sim}(\spec(\QGr A),\cc O_{\QGr A}).$$
\end{thm}

\begin{proof}
Fix an open subset $U\subseteq\spec(\Rfp)$ and consider the
composition of the functors
   $$F:\Rfp\xrightarrow{\wt{(-)}}\Qcoh\spec
     R\xrightarrow{(-)|_{f^{-1}(U)}}\Qcoh f^{-1}(U).$$
Here, for any $R$-module $M$, we denote by $\wt M$ its associated
sheaf. By definition, the stalk of $\wt M$ at a prime $P$ equals the
localized module $M_P$. We claim that $F$ annihilates $\cc S_U$. In
fact, $M\in\cc S_U$ implies $f^{-1}(\supp(M))\cap
f^{-1}(U)=\emptyset$ and therefore $\supp_R(M)\cap
f^{-1}(U)=\emptyset$. Thus $M_{P}=0$ for all $P\in f^{-1}(U)$ and
therefore $F(M)=0$. It follows that $F$ factors through $\Rfp/\cc
S_U$ and induces a map $\End_{\Rfp/\cc S_U}(R)\to\cc
O_{R}(f^{-1}(U))$ which extends to a map $\cc O_{\Rfp}(U)\to\cc
O_{R}(f^{-1}(U))$. This yields the morphism of sheaves
$f^\sharp\colon\cc O_{\Rfp}\to f_*\cc O_{R}$.

By the above $f^\sharp$ induces an isomorphism $f^\sharp_P\colon\cc
O_{\Rfp,f(P)}\to\cc O_{R,P}$ at each point $P\in\spec R$. We
conclude that $f^\sharp_P$ is an isomorphism. It follows that $f$ is
an isomorphism of ringed spaces if the map $f:\spec R\to\spec(\Rfp)$
is a homeomorphism. This last condition is a consequence of
Propositions~\ref{pr:openlattice} and~\ref{prr}. The same arguments
apply to show that
   $$f:(\Proj A,\cc O_{\Proj A})\lra{\sim}(\spec(\QGr A),\cc O_{\QGr A})$$
is an isomorphism of ringed spaces.
\end{proof}

\section{Appendix}\label{poss}

A subcategory $\cc S$ of a Grothendieck category $\cc C$ is said to
be {\em Serre\/} if for any short exact sequence
   $$0\to{X'}\to X\to {X''}\to 0$$
$X',X''\in\cc S$ \ifff $X\in\cc S$. A Serre subcategory $\cc S$ of
$\cc C$ is said to be a {\em torsion class\/} if $\cc S$ is closed
under taking coproducts. An object $C$ of $\cc C$ is said to be {\it
torsionfree\/} if $(\cc S, C)=0$.  The pair consisting of a torsion
class and the corresponding class of torsionfree objects is referred
to as a {\em torsion theory}.   Given a torsion class $\cc S$ in
$\cc C$ the {\it quotient category $\cc C/\cc S$\/} is the full
subcategory with objects those torsionfree objects $C\in\cc C$
satisfying ${\rm Ext}^1(T,C)=0$ for every $T\in \cc S$.  The
inclusion functor $i:\cc S\to\cc C$ admits the right adjoint $t:\cc
C\to\cc S$ which takes every object $X\in\cc C$ to the maximal
subobject $t(X)$ of $X$ belonging to $\cc S$. The functor $t$ we
call the {\it torsion functor}.  Moreover, the inclusion functor
$i:\cc C/\cc S\to\cc C$ has a left adjoint, the {\it localization
functor\/} ${(-)}_{\cc S}:\cc C\to\cc C/\cc S$, which is also exact.
Then,
   $$\Hom_{\cc C}(X,Y)\cong\Hom_{\cc C/\cc S}(X_{\cc S},Y)$$
for all $X\in\cc C$ and $Y\in\cc C/\cc S$. A torsion class $\cc S$
is {\it of finite type\/} if the functor $i:\cc C/\cc S\to\cc C$
preserves directed sums. If $\cc C$ is a locally coherent
Grothendieck category then $\cc S$ is of finite type \ifff $i:\cc
C/\cc S\to\cc C$ preserves direct limits (see, e.g.,~\cite{G}).

Let $\cc C$ be a Grothendieck category having a family of finitely
generated projective generators $\cc A=\{P_i\}_{i\in I}$. Let $\ff
F=\bigcup_{i\in I}\ff F^i$ be a family of subobjects, where each
$\ff F^i$ is a family of subobjects of $P_i$. We refer to $\ff F$ as
a {\em Gabriel filter\/} if the following axioms are satisfied:

\begin{itemize}
\item[$T1.$] $P_i\in\ff F^i$ for every $i\in I$;
\item[$T2.$] if $\ff a\in\ff F^i$ and $\mu:{P_j}\to{P_i}$
    then $\{\ff a:\mu\}=\mu^{-1}(\ff a)\in\ff F^j$;
\item[$T3.$] if $\ff a$ and $\ff b$ are subobjects of $P_i$ such that
    $\ff a\in\ff F^i$ and $\{\ff b:\mu\}\in\ff F^j$
    for any $\mu:{P_j}\to{P_i}$ with $\im\mu\subset\ff a$
    then $\ff b\in\ff F^i$.
\end{itemize}

\noindent In particular each $\ff F^i$ is a filter of subobjects of
$P_i$. A Gabriel filter is {\it of finite type} if each of these
filters has a  cofinal set of finitely generated objects (that is,
if for each $i$ and each $\ff a \in \ff F_i$ there is a finitely
generated $\ff b \in \ff F_i$ with $\ff a \supseteq \ff b$).

Note that if $\cc A=\{A\}$ is a ring and $\ff a$ is a right ideal of $A$, then
for every endomorphism $\mu:A\to A$
   $$\mu^{-1}(\ff a)=\{\ff a:\mu(1)\}=\{a\in A\mid \mu(1)a\in\ff a\}.$$
On the other hand, if $x\in A$, then $\{\ff a:x\}=\mu^{-1}(\ff a)$,
where $\mu\in\End A$ is such that $\mu(1)=x$.

It is well-known (see, e.g., \cite{G}) that the map
   $$\cc S\longmapsto\ff F(\cc S)=\{\ff a\subseteq P_i\mid i\in I,\, P_i/\ff a\in\cc S\}$$
establishes a bijection between the Gabriel filters (respectively
Gabriel filters of finite type) and the torsion classes on $\cc C$
(respectively torsion classes of finite type on $\cc C$).

\end{document}